\documentclass[11pt]{article}
\usepackage[left=1in,right=1in,top=1in,bottom=1in]{geometry}
\usepackage{times}
\usepackage{expl3}
\usepackage{cite}
\usepackage[table]{xcolor}
\usepackage{multirow}
\usepackage{stackengine} 
\usepackage{hhline}
\usepackage{lipsum}
\usepackage{titlesec}
\usepackage{wrapfig}
\usepackage{epsfig}
\usepackage{graphicx}
\usepackage{amsmath}
\usepackage[title]{appendix}
\usepackage{amssymb}
\usepackage{epstopdf}
\usepackage{boldline}
\usepackage{calligra}
\usepackage{url}
\usepackage{blindtext}

\newcommand{\define}{\stackrel{\mbox{\tiny def}}{=}}

\newtheorem{definition}{Definition}
\newtheorem{theorem}{Theorem}

\newtheorem{lemma}{Lemma}

\usepackage{mathtools}
\usepackage{epstopdf}
\usepackage{balance}
\usepackage{thmtools}
\usepackage{thm-restate}
\usepackage{hyperref}
\usepackage{cleveref}

\usepackage[ruled,vlined]{algorithm2e}
\include{pythonlisting}
\newcommand{\ostar}{\mathbin{\mathpalette\make@circled\star}}

\makeatletter
\newcommand{\removelatexerror}{\let\@latex@error\@gobble}
\makeatother
\makeatletter
\newcommand*{\rom}[1]{\expandafter\@slowromancap\romannumeral #1@}
\makeatother

\ExplSyntaxOn
\newcommand\latinabbrev[1]{
  \peek_meaning:NTF . {
    #1\@}%
  { \peek_catcode:NTF a {
      #1.\@ }%
    {#1.\@}}}
\ExplSyntaxOff

\titleclass{\subsubsubsection}{straight}[\subsubsection]
\date{}
\begin{document}
\vspace{1cm}
\title{Geometric Approach For Majorizing Measures}\vspace{1.8cm}
\author{Shih~Yu~Chang 
\thanks{Shih Yu Chang is with the Department of Applied Data Science,
San Jose State University, San Jose, CA, U. S. A. (e-mail: {\tt
shihyu.chang@sjsu.edu}).
           }}

\maketitle

\begin{abstract}
Gaussian processes can be considered as subsets of a standard Hilbert space, but the geometric understanding that would relate the size of a set with the size of its convex hull is still lacking. In this work, we adopt a geometric approach to the majorizing measure problem by identifying the covering number relationships between a given space $T$ and its convex hull, represented by $T_h$. If the space $T$ is a closed bounded polyhedra in $\mathbb{R}^n$, we can evaluate the volume ratio between the space $T$ and its convex hull obtained by the Quickhull algorithm. If the space $T$ is a general compact object in $\mathbb{R}^n$ with non-empty interior, we first establish a more general reverse Brunn-Minkowski inequality for nonconvex spaces which will assist us to bound the volume of $T_h$ in terms of the volume of $T$ if $T_h$ can be acquired by the finite average of the space $T$ with respect to the Minkowski sum. If the volume ratio between the space $T$ and the space $T_h$ is obtained, the covering number ratio between the space $T$ and the space $T_h$ can also be obtained which will be used to build majorizing measure inequality. For infinite dimensional space, we show that the constant $L$ at majorizing measure inequality may not always exist and the existence condition will depend on geometric properties of $T$. 
\end{abstract}


\section{Introduction}\label{sec:Introduction}

In order to estimate the supremum of a stochastic process, generic chaining method
is a powerful tool for this aspect. This method stemmed from the classical chaining method and the later majorizing measures method, which were developed by several great scholars, Kolmogorov, Dudley, Fernique and Talagrand,
to understand the supremum properties of stochastic processes. 

Let us recall the following definition.
\begin{definition}\label{def:gamma 2}
Given a metric space $(T, d(\cdot, \cdot))$, we define $\gamma_{\alpha}(T,d(\cdot, \cdot))$ as  
\begin{eqnarray}
\gamma_\alpha(T, d(\cdot, \cdot)) =\inf \sup\limits_{t \in T} \sum\limits_{m \geq 0}2^{m/\alpha} \Delta(\mathrm{A}_n(t)),  
\end{eqnarray}
where the infimum is taken over all admissible sequences and $\Delta(\mathrm{A}_n(t))$ denotes the diameter of $\mathrm{A}_n(t)$ with respect to the metric $d(\cdot, \cdot)$. 
\end{definition}
In the Gaussian case, that is when $X_t = \sum\limits_{t_k \in T} t_k g_k$, where $g_k$ are i.i.d. standard Gaussians, we have the celebrated Fernique–Talagrand majorizing measure theorem as:
\begin{eqnarray}
\frac{1}{L} \gamma_2(T, d(\cdot, \cdot)) \leq \mathbb{E} \sup\limits_{t \in T} X_t   \leq L \gamma_2(T, d(\cdot, \cdot)).
\end{eqnarray}
The majorizing measure theorem is a main method used to prove Theorem 2.11.1 in~\cite{talagrand2021upper}. 
Gaussian processes can be seen as subsets of a standard Hilbert space, but the
geometric understanding that would relate the size of a set with the size of its
convex hull is still insufficient~\cite{talagrand2021upper}. The purpose of this work is 
to establish Theorem 2.11.1 in~\cite{talagrand2021upper} based on a geometrical approach, i.e., we wish to show
\begin{eqnarray}\label{eq:2.130}
\gamma_2 (T_h) \leq L \gamma_2(T),
\end{eqnarray}
where $T_h$ is the convex hull for the original space $T$.

We consider both finite dimensional and infinite dimensional situations. For a finite dimensional situation, we consider two subcases: polyhedra and general object in $\mathbb{R}^n$. If the space $T$ is a closed bounded polyhedra in $\mathbb{R}^n$, we determine the volume ratio between the space $T$ and its convex hull obtained by the Quickhull algorithm, denoted as $T_h$. If the space $T$ is a general compact object in $\mathbb{R}^n$ with non-empty interior, we first prepare a more general reverse Brunn-Minkowski inequality for nonconvex spaces which will help us to bound the volume of $T_h$ in terms of the volume of $T$ if $T_h$ can be obtained by the finite average of the space $T$ with respect to the Minkowski sum. If the volume ratio between the space $T$ and the space $T_h$ is obtained, the covering number ratio between the space $T$ and the space $T_h$ can also be obtained which will be used to establish Theorem 2.11.1. For infinite dimensional space, since the important part in the geometric approach 
to prove Theorem 2.11.1 in~\cite{talagrand2021upper} is to know the the covering number (metric entropy) relations between the space $T$ and its convex hull space $T_h$, the authors in~\cite{cockreham2017metric} have results about the estimation of the covering number of $T_h$ given the covering number of $T$. We demonstrate that the constant $L$ at Theorem 2.11.1 in~\cite{talagrand2021upper} may not always exist and the existence will depend on geometric properties of $T$. 
   
The rest of this paper is organized as follows. The volume ratio for the polyhedra and the general object in $\mathbb{R}^n$ is discussed in Section~\ref{sec:Volume Difference Between T and Th}. In Section~\ref{sec:Covering Number Bounds by Volume}, we bound the covering number of $T_h$ in terms of the covering number of $T$ for polyhedra and general situations. Majorizing measuring theorem is proved from a geometric perspective in Section~\ref{sec:Majorizing Measures By Geometry Approach}. For the space with infinite dimension, we provide the condition for the existence of the constant $L$ in Eq.~\eqref{eq:2.130} at Section~\ref{sec:Infinite Dimensional T Discussions}. 

\section{Volume Ratio Between $T$ and $T_h$}\label{sec:Volume Difference Between T and Th}
 
In this section, we will consider the volume ratio between the original geometric space and its convex hull. All objects are assumed to be compact with non-empty interior. In Section~\ref{sec:Polyhedra Case}, we will consider $T$ as a polyhedra object and the general geometric objects in Eucleduean space is discussed in Section~\ref{sec:General Case}. If the space $T$ is a convex space, its convex hull space, denoted as $T_h$, is idential to $T$. Therefore, we will focus the volume difference for the situation that $T$ is not a convex space. 

\subsection{Polyhedra Case}\label{sec:Polyhedra Case}

We assume that the space $T$ is a closed bounded polyhedra in $\mathbb{R}^n$ having an orientable boundary, denoted as $\partial T$, which can be triangulated into a set of $(n-1)$-dimensional simplices. From Eq. (2.1) in~\cite{allgower1986computing}, we have the volume of $T$ expressed as
\begin{eqnarray}\label{eq1:volume formula 1}
\mbox{Vol}_n\left( T \right) &=& \sum\limits_{\sigma \in \partial T} \frac{1}{n!} \mbox{det}\left(v_1(\sigma), \cdots, v_n(\sigma) \right),
\end{eqnarray}
where $v_i(\sigma)$ are the vertices of the $(n-1)$-simplex $\sigma$ ordered according to the orientation of $\sigma$. The other volume formula provided by Eq. (2.2) in~\cite{allgower1986computing} can be expressed as 
\begin{eqnarray}\label{eq1:volume formula 2}
\mbox{Vol}_n\left( T \right) &=& (-1)^{(n-1}\sum\limits_{\sigma \in \partial T} \left( \frac{1}{n}\sum\limits_{i=1}^{n} v_i[n]  \right) \frac{1}{(n-1)!} \mbox{det}\left(     \begin{array}{ccc}
       1  & \cdots & 1   \\
       v_1[\neg n]  & \cdots & v_n[\neg n]  \\
    \end{array} \right),
\end{eqnarray}
where $v_i[n]$ is the $n$-th coordinate of the point $v_i$, and $v_i[\neg n]$ is obtained by deleting the $n$-th coordinate from $v_i$.  

Following Lemma~\ref{lma:vol ratio polyhedra} is provided to give the volume ratio between $T_h$ and $T$.
\begin{lemma}\label{lma:vol ratio polyhedra}
Given a compact and connected polyhedra $T$ in $\mathbb{R}^n$, we have volume ratio between $T_h$ and $T$, denoted as $R_{\mbox{Poly},n}$, which can be expressed as
\begin{eqnarray}
R_{\mbox{Poly},n} = \frac{ \mbox{Vol}_n\left( T_h \right)  }{  \mbox{Vol}_n\left( T \right)  },
\end{eqnarray}
where $\mbox{Vol}_n\left( T_h \right)$ and $\mbox{Vol}_n\left( T \right)$ are volumes of $T_h$ and $T$.
\end{lemma}

\textbf{Proof:}
We can apply the Quickhull algorithm to convert $T$ into $T_h$, see~\cite{barber1996quickhull}. From Theorem 2.3 in~\cite{barber1996quickhull}, they show that the Quickhull algorithm can produce the convex hull of a set of points (polyhedra) in $\mathbb{R}^n$. Then, we can apply Eq.~\eqref{eq1:volume formula 1} or Eq.~\eqref{eq1:volume formula 2} to evaluate $\mbox{Vol}_n\left( T_h \right) $. Then, this Lemma is proved. 
$\hfill \Box$

Note that $R_{\mbox{Poly},n} \geq 1$ and the equality holds when $T$ is a convex set.

\subsection{General Case}\label{sec:General Case}

In this section, we will extend polyhedra to general objects in $\mathbb{R}^n$. Our main tool is to reverse Brunn–Minkowski theorem. Let us recall the traditional reverse Brunn–Minkowski inequality first. There is a constant $C_1$, independent of $n$, such that for any two centrally symmetric convex objects $A$ and $B$ in $\mathbb{R}^n$, there are volume-preserving linear maps $\rho$ and $\varrho$ from $\mathbb{R}^{n}$ to $\mathbb{R}^{n}$ such that for any real numbers $s,t >0$ and any natural number $m$, we have
\begin{eqnarray}\label{eq:old rev-BM inequality}
[\mbox{Vol}_n \left( s \rho(A) \oplus t \varrho(B) \right)]^{1/m} \leq C_1 \left( s [\mbox{Vol}_n ( \rho(A) )]^{1/m} +  t [\mbox{Vol}_n ( \varrho(B) )]^{1/m}   \right),
\end{eqnarray}
where $\oplus$ is the Minkowski sum operator between two geometric objects. However, traditional reverse Brunn–Minkowski theorem requires that the objects $A$ and $B$ are centrally symmetric convex bodies. We will prove a more general Theorem~\ref{thm:rev-BM Inequality} that relaxs the objects $A$ and $B$ required to be  centrally symmetric convex objects.

Given a compact object $D$ in $\mathbb{R}^n$, we have its circumscribed ball, denoted as $S_{D}$, and the ratio between $\mbox{Vol}(S_{D})$ and $\mbox{Vol}(D)$ is represented by $\beta_D$, i.e., 
\begin{eqnarray}\label{eq:gamma vol ratio circum}
\beta_D = \frac{\mbox{Vol}(S_{D})}{\mbox{Vol}(D)}.
\end{eqnarray}
Then, we have the following theorem about more general class of objects in $\mathbb{R}^n$.
\begin{theorem}\label{thm:rev-BM Inequality}
Given two compact objects with non-empty interior in $\mathbb{R}^n$, denoted as $A$ and $B$, two volume-preserving linear maps $\rho, \varrho$, any real numbers $s,t >0$ and any natural number $m$,  we have
\begin{eqnarray}\label{eq1:thm:rev-BM Inequality}
[\mbox{Vol}_n \left( s \rho(A) \oplus t \varrho(B) \right)]^{1/m} \leq C_1 \left( s [ \beta_{\rho(A)} \mbox{Vol}_n ( \rho(A) )]^{1/m} +  t [\beta_{\varrho(B)}\mbox{Vol}_n ( \varrho(B) )]^{1/m}   \right),
\end{eqnarray}
where $C_2$ is a constant.
\end{theorem}

\textbf{Proof:}
Since we have
\begin{eqnarray}\label{eq2:thm:rev-BM Inequality}
[\mbox{Vol}_n \left( s \rho(A) \oplus t \varrho(B) \right)]^{1/m} &\leq_1& 
[\mbox{Vol}_n \left( s S_{\rho(A)} \oplus t S_{\varrho(B)} \right)]^{1/m} \nonumber \\
&\leq_2 &  C_1 \left( s [\mbox{Vol}_n (  S_{\rho(A)} )]^{1/m} +  
t [\mbox{Vol}_n (   S_{\varrho(B)} )]^{1/m}   \right) \nonumber \\
& =  &  C_1 \left( s [ \beta_{\rho(A)} \mbox{Vol}_n ( \rho(A) ) ]^{1/m} +  
t [ \beta_{\varrho(B)} \mbox{Vol}_n ( \varrho(B) ) ]^{1/m}   \right),
\end{eqnarray}
where the inequality $\leq_1$~\footnote{We can consider to apply other center symmetric shapes, like $p$-ball, to provide better circumscribed fit of the original shapes $\rho(A)$ or $\varrho(B)$ in $\leq_1$ of Eq.~\eqref{eq2:thm:rev-BM Inequality}, see~\cite{bastero1995extension, litvak34constant}. This means that we can provide tighter bound in $\leq_1$ of Eq.~\eqref{eq2:thm:rev-BM Inequality} according to the geometric properties of $\rho(A)$ or $\varrho(B)$.} comes from the fact that the volume of an object is less or equal than its circumscribed ball;  the inequality $\leq_2$ comes from Eq.~\eqref{eq:old rev-BM inequality}; the lasy equality comes from the definition $\beta$ provided by Eq.~\eqref{eq:gamma vol ratio circum}.
$\hfill \Box$

Given an object $A$ in $\mathbb{R}^n$, we define $A(k)$ as 
\begin{eqnarray}
A(k) = \frac{1}{k}\oplus_{i=1}^k A_i,
\end{eqnarray} 
where $A_i = A$. Then, we have the following Lemma about the volume $A(k)$ and the volume $A(k-1)$.

\begin{lemma}\label{lma:vol A k relation}
Let $A$ be a compact objects with non-empty interior in $\mathbb{R}^n$, then we have 
\begin{eqnarray}\label{eq1:lma:vol A k relation}
\mbox{Vol}_n(A(k)) \leq \left( \frac{2C^{(k-1)}_2}{k} + \frac{C_2 (  C^{(k-2)}_2  - 1 )   }{k (C_2 - 1)  }  \right)
\mbox{Vol}_n(A),
\end{eqnarray}
where $C_2 = C_1 C'_2$ and $C'_2$ is an upper bound for $\beta_{A(k)}$. 
\end{lemma}

\textbf{Proof:}
From~\cite{fradelizi2018convexification}, we know that $A(k) \rightarrow A_h$ when $k \rightarrow \infty$, where $A_h$ is the convex hull of $A$. Since $A_h$ is also a volume bounded object, we can upper bound the ratio parameter $\beta_{A(k)}$ by $C'_2$ for all $k$, where $C'_2$ is a positive number. 

Because we have 
\begin{eqnarray}\label{eq2:lma:vol A k relation}
\mbox{Vol}_n(A(k)) &=& \mbox{Vol}_n \left( \frac{k-1}{k}A(k-1) + \frac{1}{k} A \right) \nonumber \\
&\leq_1& C_1 \left(  \frac{k-1}{k} [ \beta_{A(k-1) } \mbox{Vol}_n ( A(k-1) )] +  \frac{1}{k} [\beta_{A}\mbox{Vol}_n ( A )]  \right) \nonumber \\
&\leq_2& C_1 \left(  \frac{k-1}{k} [ C'_2 \mbox{Vol}_n ( A(k-1) )] +  \frac{1}{k} [C'_2 \mbox{Vol}_n ( A )]  \right) 
\end{eqnarray}
where the inequality $\leq_1$ is obtained from Theorem~\ref{thm:rev-BM Inequality} by setting $m=1$; the inequality $\leq_2$ is obtained by using the upper bound $C'_2$ for $\beta_{A(k-1)}$ and $\beta_{A}=\beta_{A(1)}$. 

By applying Eq.~\eqref{eq2:lma:vol A k relation} iterativelys, we can uppder the volume of $A(k)$ by the volume of $A$ as 
\begin{eqnarray}\label{eq3:lma:vol A k relation}
\mbox{Vol}_n(A(k)) \leq \left( \frac{2C^{(k-1)}_2}{k} + \frac{C_2 (  C^{(k-2)}_2  - 1 )   }{k (C_2 - 1)  }  \right)
\mbox{Vol}_n(A), 
\end{eqnarray}
where $C_2 = C_1 C'_2$. This Lemma is proved. 
$\hfill \Box$

Following Lemma~\ref{lma:vol ratio general} is provided to give the volume ratio between $T_h$ and $T$ when $T$ is a general compact object with non-empty interior in $\mathbb{R}^n$. 
\begin{lemma}\label{lma:vol ratio general}
Given $T$ as a compact object with non-empty interior in $\mathbb{R}^n$ and we assume that $A(k_h) = A_h$  for some natural number $k_h$. This means that $A(k_h+n)$ will be equal to $A_h$ for any natural number $n$ since $A(k_h)$ is a convex set already. We have volume ratio between $T_h$ and $T$, denoted as $R_{\mbox{Gen},n}$, which can be expressed as
\begin{eqnarray}
R_{\mbox{Gen},n} &=& \frac{ \mbox{Vol}_n\left( T_h \right)  }{  \mbox{Vol}_n\left( T \right)  } \nonumber \\
&\leq& \left( \frac{2C^{(k_h-1)}_2}{k_h} + \frac{C_2 (  C^{(k_h-2)}_2  - 1 )   }{k_h (C_2 - 1)  }  \right).
\end{eqnarray}
\end{lemma}
\textbf{Proof:}
This Lemma can be obtained directly from Lemma~\ref{lma:vol A k relation}.  
$\hfill \Box$

\section{Covering Number Bounds by Volume}\label{sec:Covering Number Bounds by Volume}

In this section, we will try to upper bound the covering number of $T_h$ by the covering number of $T$ via covering number bounds by volume. We begin with the covering number definition. 
\begin{definition}
We use $N(A, d(\cdot, \cdot), \epsilon)$ to represent the covering number of the space $A$ with $\epsilon$-ball with distance metric function $d(\cdot, \cdot)$. Then, $N(A, d(\cdot, \cdot), \epsilon)$ can be defined as
\begin{eqnarray}
N(T, d(\cdot, \cdot), \epsilon) \define \min\{n: \mbox{the number of $\epsilon$-ball to cover $T$} \}.
\end{eqnarray}
\end{definition}

From Theorem 14.2 in~\cite{CoveringAndPacking2016}, we have
\begin{eqnarray}\label{eq:CoveringAndPacking2016}
\left(\frac{1}{\epsilon}\right)^n \frac{\mbox{Vol}(A) }{ \mbox{Vol}(\mathbf{B}) } \leq N(A, d(\cdot, \cdot), \epsilon) \leq \frac{\mbox{Vol}(A \oplus \frac{\epsilon}{2} \mathbf{B} ) }{ \mbox{Vol}( \frac{\epsilon}{2} \mathbf{B}) }
\leq_1 \left(\frac{3}{\epsilon}\right)^n \frac{\mbox{Vol}(A) }{ \mbox{Vol}(\mathbf{B}) }
\end{eqnarray}
where $\mathbf{B}$ is the unit norm ball, and the inequality $\leq_1$ is valid if $A$ is a convex set and $\epsilon \mathbf{B} \subset A$. 

We will have the following two Lemmas to upper bound the covering number of $T_h$ by the covering number of $T$ with respect to different types of $T$. 

\begin{lemma}\label{lma:Th convering num upper bound poly}
Given a compact and connected polyhedra $T$ in $\mathbb{R}^n$, we have the following relation between the covering number of $T_h$ and the covering number of $T$:
\begin{eqnarray}\label{eq1:lma:Th convering num upper bound poly}
N(T_h, d(\cdot, \cdot), \epsilon) &\leq&  R_{\mbox{Poly},n} 3^n  N(T, d(\cdot, \cdot), \epsilon),
\end{eqnarray}
where $R_{\mbox{Poly},n}$ is the volume ratio given by Lemma~\ref{lma:vol ratio polyhedra}.
\end{lemma}
\textbf{Proof:}

Because $T_h$ is a convex set, we have
\begin{eqnarray}
N(T_h, d(\cdot, \cdot), \epsilon) &\leq&   \left(\frac{3}{\epsilon}\right)^n \frac{\mbox{Vol}(T_h) }{ \mbox{Vol}(\mathbf{B}) } \nonumber \\
&=_1&   \left(\frac{3}{\epsilon}\right)^n \frac{R_{\mbox{Poly},n} \mbox{Vol}(T) }{ \mbox{Vol}(\mathbf{B}) } \nonumber \\
&=&  R_{\mbox{Poly},n} 3^n  \left(\frac{1}{\epsilon}\right)^n \frac{ \mbox{Vol}(T) }{ \mbox{Vol}(\mathbf{B}) } \nonumber \\
&\leq & R_{\mbox{Poly},n} 3^n   N(T, d(\cdot, \cdot), \epsilon),
\end{eqnarray}
where the equality $=_1$ comes from Lemma~\ref{lma:vol ratio polyhedra}, and the first and the last 
inequalities are obtained from Eq.~\eqref{eq:CoveringAndPacking2016}. 
$\hfill \Box$

\begin{lemma}\label{lma:Th convering num upper bound gen}
Given a compact object $T$ with non-empty interior in $\mathbb{R}^n$ and assume that $T_h = \frac{1}{k}\oplus_{i=1}^{k} T$ with finite $k$, we have the following relation between the covering number of $T_h$ and the covering number of $T$:
\begin{eqnarray}\label{eq1:lma:Th convering num upper bound gen}
N(T_h, d(\cdot, \cdot), \epsilon) &\leq&  R_{\mbox{Gen},n} 3^n  N(T, d(\cdot, \cdot), \epsilon),
\end{eqnarray}
where $R_{\mbox{Gen},n}$ is the volume ratio given by Lemma~\ref{lma:vol ratio general}.
\end{lemma}
\textbf{Proof:}

Because $T_h$ is a convex set, we have
\begin{eqnarray}
N(T_h, d(\cdot, \cdot), \epsilon) &\leq&   \left(\frac{3}{\epsilon}\right)^n \frac{\mbox{Vol}(T_h) }{ \mbox{Vol}(\mathbf{B}) } \nonumber \\
&=_1&   \left(\frac{3}{\epsilon}\right)^n \frac{R_{\mbox{Gen},n} \mbox{Vol}(T) }{ \mbox{Vol}(\mathbf{B}) } \nonumber \\
&=&  R_{\mbox{Gen},n} 3^n  \left(\frac{1}{\epsilon}\right)^n \frac{ \mbox{Vol}(T) }{ \mbox{Vol}(\mathbf{B}) } \nonumber \\
&\leq & R_{\mbox{Gen},n} 3^n   N(T, d(\cdot, \cdot), \epsilon),
\end{eqnarray}
where the equality $=_1$ comes from Lemma~\ref{lma:vol ratio general}, and the first and the last 
inequalities are obtained from Eq.~\eqref{eq:CoveringAndPacking2016}. 
$\hfill \Box$

\section{Majorizing Measures By Geometry Approach}\label{sec:Majorizing Measures By Geometry Approach}

From Lemma~\ref{lma:Th convering num upper bound poly} and Lemma~\ref{lma:Th convering num upper bound gen}, we are ready to prove Theorem 2.11.1 in~\cite{talagrand2021upper} geometrically at Eucledean space with dimension $n$. 

\begin{theorem}\label{thm:eq 2.130 in talagrand2021upper Book Poly}
Given the space $T$ as a closed bounded polyhedra in $\mathbb{R}^n$ having an orientable boundary, we have 
\begin{eqnarray}
\gamma_\alpha(T_h, d(\cdot, \cdot)) \leq L_{\mbox{Poly}} \gamma_\alpha(T, d(\cdot, \cdot)),
\end{eqnarray}
where $L_{\mbox{Poly}}$ is the constant depending on the underlying geometry of the space $T$ and $\alpha$. 
\end{theorem}
\textbf{Proof:}

Because we have 
\begin{eqnarray}
\gamma_{\alpha}(T_h, d(\cdot, \cdot)) & \leqslant_{\alpha} & \int_0^{\infty}   \left(\log   N(T_h, d(\cdot, \cdot), \epsilon) \right)^{1/\alpha}   d \epsilon \nonumber \\
&\leq_1 &  \int_0^{\infty}   \left(\log  R_{\mbox{Poly},n} 3^n N(T, d(\cdot, \cdot), \epsilon) \right)^{1/\alpha}   d \epsilon \nonumber \\
&\leq &  \left(\frac{\log (R_{\mbox{Poly},n} 3^n)}{\log 2 }  + 1 \right)^{1/\alpha}\int_0^{\infty}   \left(\log  N(T, d(\cdot, \cdot), \epsilon) \right)^{1/\alpha}   d \epsilon \nonumber \\
& \leqslant_{\alpha}& \left(\frac{\log (R_{\mbox{Poly},n} 3^n)}{\log 2 }  + 1 \right)^{1/\alpha}\gamma_{\alpha}(T, d(\cdot, \cdot)), 
\end{eqnarray}
where the inequality $\leq_1$ comes from Lemma~\ref{lma:Th convering num upper bound poly}, and the first and last inequalities come from Theorem 1.2 in~\cite{talagrand2001majorizing}. This Lemma is proved. 
$\hfill \Box$

\begin{theorem}\label{thm:eq 2.130 in talagrand2021upper Book Gen}
Given a compact object $T$ in $\mathbb{R}^n$ such that $T_h = \frac{1}{k}\oplus_{i=1}^k T$ with finite $k$, we have 
\begin{eqnarray}
\gamma_\alpha(T_h, d(\cdot, \cdot)) \leq L_{\mbox{Gen}} \gamma_\alpha(T, d(\cdot, \cdot)),
\end{eqnarray}
where $L_{\mbox{Gen}}$ is the constant depending on the underlying geometry of the space $T$ and $\alpha$. 
\end{theorem}
\textbf{Proof:}

Since we have 
\begin{eqnarray}
\gamma_{\alpha}(T_h, d(\cdot, \cdot)) & \leqslant_{\alpha} & \int_0^{\infty}   \left(\log   N(T_h, d(\cdot, \cdot), \epsilon) \right)^{1/\alpha}   d \epsilon \nonumber \\
&\leq_1 &  \int_0^{\infty}   \left(\log  R_{\mbox{Gen},n} 3^n N(T, d(\cdot, \cdot), \epsilon) \right)^{1/\alpha}   d \epsilon \nonumber \\
&\leq &  \left(\frac{\log (R_{\mbox{Gen},n} 3^n)}{\log 2 }  + 1 \right)^{1/\alpha}\int_0^{\infty}   \left(\log  N(T, d(\cdot, \cdot), \epsilon) \right)^{1/\alpha}   d \epsilon \nonumber \\
& \leqslant_{\alpha}& \left(\frac{\log (R_{\mbox{Gen},n} 3^n)}{\log 2 }  + 1 \right)^{1/\alpha}\gamma_{\alpha}(T, d(\cdot, \cdot)), 
\end{eqnarray}
where the inequality $\leq_1$ comes from Lemma~\ref{lma:Th convering num upper bound gen}, and the first and last inequalities come from Theorem 1.2 in~\cite{talagrand2001majorizing}. This Lemma is proved. 
$\hfill \Box$

\section{Infinite Dimensional $T$ Discussion}\label{sec:Infinite Dimensional T Discussions}

From Section~\ref{sec:Majorizing Measures By Geometry Approach}, the key part in the geometric approach 
to prove Theorem 2.11.1 in~\cite{talagrand2021upper} is to represent the covering number $N(T_h, d(\cdot, \cdot), \epsilon)$ in terms of the covering number $N(T, d(\cdot, \cdot), \epsilon)$. We will show that, when $T$ is a subset of an infinite dimensional space, e.g., Banach space or Hilbert space, the constant $L$ at Eq.~\eqref{eq:2.130} may not always exist. The existence will depend on geometric properties of $T$.

According to Section 2 in~\cite{cockreham2017metric}, the authors list some specific estimates
of $\log N(T_h, d(\cdot, \cdot), \epsilon)$ under various common assumptions on $\log N(T, d(\cdot, \cdot), \epsilon)$. Here, we assume $p=2$ and $q=1$ (convex hull) in those parameters used in Section 2 in~\cite{cockreham2017metric}.

If $\log N(T, d(\cdot, \cdot), \epsilon) = O(\epsilon^{- \chi} \left\vert \log \epsilon \right\vert^{\psi}   )$ for $\chi > 2$, and $\psi \in \mathbb{R}$, we have $\log N(T_h, d(\cdot, \cdot), \epsilon) = O(\epsilon^{- \chi} \left\vert \log \epsilon \right\vert^{\psi}  )$. In this situation, we can find a constant, say $C_3$, to upper bound the following ratio between the  covering number $N(T_h, d(\cdot, \cdot), \epsilon)$ and the covering number $N(T, d(\cdot, \cdot), \epsilon)$:
\begin{eqnarray}\label{eq:case 1}
\frac{\log N(T_h, d(\cdot, \cdot), \epsilon) }{ \log N(T, d(\cdot, \cdot), \epsilon) } \leq C_3.
\end{eqnarray}
According to the proof in Theorem~\ref{thm:eq 2.130 in talagrand2021upper Book Poly} (or Theorem~\ref{thm:eq 2.130 in talagrand2021upper Book Gen}), the ratio between $ \frac{ \gamma_\alpha(T_h, d(\cdot, \cdot)) }{ \gamma_\alpha(T, d(\cdot, \cdot)) }$  can be upper bound by the integration with respect to $\epsilon$ from 0 to $\Delta(T)$ for contant integrand $C_3$. 

If $\log N(T, d(\cdot, \cdot), \epsilon) = O(\epsilon^{- 2} \left\vert \log \epsilon \right\vert^{\psi}   )$ for $\psi > -2$, we have $\log N(T_h, d(\cdot, \cdot), \epsilon) = O(\epsilon^{- 2} \left\vert \log \epsilon \right\vert^{2 + \psi}  )$. In this situation, we can find a function of $\epsilon$ to upper bound the following ratio between the  covering number $N(T_h, d(\cdot, \cdot), \epsilon)$ and the covering number $N(T, d(\cdot, \cdot), \epsilon)$:
\begin{eqnarray}\label{eq:case 2}
\frac{\log N(T_h, d(\cdot, \cdot), \epsilon) }{ \log N(T, d(\cdot, \cdot), \epsilon) } \leq  C_4 \left\vert \log \epsilon \right\vert^{2}.
\end{eqnarray}
According to the proof in Theorem~\ref{thm:eq 2.130 in talagrand2021upper Book Poly} (or Theorem~\ref{thm:eq 2.130 in talagrand2021upper Book Gen}), the ratio between $ \frac{ \gamma_\alpha(T_h, d(\cdot, \cdot)) }{ \gamma_\alpha(T, d(\cdot, \cdot)) }$  can be upper bound by the integration with respect to $\epsilon$ from 0 to $\Delta(T)$ for the integrand as $ C_4 \left\vert \log \epsilon \right\vert^{2}$ since such integration exists.

If $\log N(T, d(\cdot, \cdot), \epsilon) = O(\epsilon^{- 2} \left\vert \log \epsilon \right\vert^{\psi}   )$ for $\psi = -3$, we have $\log N(T_h, d(\cdot, \cdot), \epsilon)=$ \\  $O(\epsilon^{- 2} \left(\left\vert \log \left\vert \log \epsilon \right\vert \right\vert \right)^{2 + \psi}  )$. In this situation, we can find a function of $\epsilon$ to upper bound the following ratio between the  covering number $N(T_h, d(\cdot, \cdot), \epsilon)$ and the covering number $N(T, d(\cdot, \cdot), \epsilon)$:
\begin{eqnarray}\label{eq:case 4}
\frac{\log N(T_h, d(\cdot, \cdot), \epsilon) }{ \log N(T, d(\cdot, \cdot), \epsilon) } \leq  C_5 \frac{\left\vert \log \epsilon \right\vert^3}{\left\vert  \log \left\vert \log \epsilon \right\vert \right\vert}.
\end{eqnarray}
According to the proof in Theorem~\ref{thm:eq 2.130 in talagrand2021upper Book Poly} (or Theorem~\ref{thm:eq 2.130 in talagrand2021upper Book Gen}), the ratio between $ \frac{ \gamma_\alpha(T_h, d(\cdot, \cdot)) }{ \gamma_\alpha(T, d(\cdot, \cdot)) }$  can be upper bound by the integration with respect to $\epsilon$ from 0 to $\Delta(T)$ for the integrand as $ C_5 \frac{\left\vert \log \epsilon \right\vert^3}{\log \left\vert \log \epsilon \right\vert}$. However, if the diameter $\Delta(T)$ equals one, this integral does not exist. 

Therefore, the existence of the constant $L$ at Theorem 2.11.1 in~\cite{talagrand2021upper} can be summarized by the following Theorem~\ref{thm:the existence of L}
\begin{theorem}\label{thm:the existence of L}
Suppose we have a function $f(\epsilon)$ which satisfies the following:
\begin{eqnarray}
\frac{\log N(T_h, d(\cdot, \cdot), \epsilon) }{ \log N(T, d(\cdot, \cdot), \epsilon) }  \leq f(\epsilon),
\end{eqnarray}
given $0 < \epsilon \leq \Delta(T)$. If the following integral exists:
\begin{eqnarray}
\int_0^{\Delta(T)} f(\epsilon),
\end{eqnarray}
we can find the constant $L$ at Eq.~\eqref{eq:2.130}. 
\end{theorem}
\textbf{Proof:}

The proof is the same as the argument used in the proof for Theorem~\ref{thm:eq 2.130 in talagrand2021upper Book Poly}
or Theorem~\ref{thm:eq 2.130 in talagrand2021upper Book Gen}.
$\hfill \Box$

\bibliographystyle{IEEETran}
\bibliography{GeometricApproach_GC_Bib}

\end{document}